# The Arrow of Time in Multivariate Time Series

**Stefan Bauer**     STEFAN.BAUER@INF.ETHZ.CH
*Department of Computer Science, ETH Zurich*
*8092 Zurich, Switzerland*

**Bernhard Schölkopf**     BS@TUEBINGEN.MPG.DE
*Max Planck Institute for Intelligent Systems*
*72076 Tübingen, Germany*

**Jonas Peters**     JONAS.PETERS@TUEBINGEN.MPG.DE
*Max Planck Institute for Intelligent Systems*
*72076 Tübingen, Germany*

## Abstract

We prove that a time series satisfying a (linear) multivariate autoregressive moving average (VARMA) model satisfies the same model assumption in the reversed time direction, too, if all innovations are normally distributed. This reversibility breaks down if the innovations are non-Gaussian. This means that under the assumption of a VARMA process with non-Gaussian noise, the arrow of time becomes detectable. Our work thereby provides a theoretic justification of an algorithm that has been used for inferring the direction of video snippets. We present a slightly modified practical algorithm that estimates the time direction for a given sample and prove its consistency. We further investigate how the performance of the algorithm depends on sample size, number of dimensions of the time series and the order of the process. An application to real world data from economics shows that considering multivariate processes instead of univariate processes can be beneficial for estimating the time direction. Our result extends earlier work on univariate time series. It relates to the concept of causal inference, where recent methods exploit non-Gaussianity of the error terms for causal structure learning.

## 1. Introduction

The goal of this work is to infer the direction of a given multivariate time series. Figure 1 shows an example of a time series in forward and backward direction. The task is to decide, which of both direction corresponds to the correct time direction. This question is mainly academic but has received a lot of attention in literature, especially in physics, see Reichenbach (1956) or Price (1996). Hernández-Lobato et al. (2011), Morales-Mombiela et al. (2013) and Hernández-Lobato et al. (2014) discuss the derivation and application of Gaussianity measures to detect the direction of causal time series. However, their approach is based on the empirical observation that the residuals of a linear fit in the forward direction are less Gaussian than the residuals in the backward direction. We are not aware of any theoretical result that clarifies under which assumption this heuristic holds. The main thrust of the present paper is to provide such identifiability results. Pickup et al. (2014) estimate the direction of a video snippet that is played either forwards or backwards. While the authors use an algorithm that is similar to the one we present in here, they do not provide a theoretical justification for their approach and refer to the univariate version of this problem (Peters et al., 2009). Our work provides a post hoc justification of Pickup et al. (2014) by proving identifiability statements. It further provides a consistency result for a version of the algorithm first derived in Peters et al. (2009); consistency was unknown both in the univariate and the multivariate case. The question of time-reversibility is furthermore interesting from a causal point of view. In causal discovery, we observe an i.i.d. sample from a multivariate distribution and try to identify the underlying causal structure. Because causes are widely accepted to precede their effects (if we have accurate and instantaneous measurements), similar principles can be applied to the problem of time reversibility. More concretely, our results are based on the Darmois-Skitovich theorem (Ibragimov, 2014), the univariate version of which is at the basis of the causal discovery method LiNGAM (Shimizu et al., 2006). Both approaches make use of non-Gaussian residuals in linear models. While it is difficult to find data sets with a known causal structure, we have access to the correct time direction of almost all observed time series. Shimizu et al. (2006) proposed to apply LiNGAM to the time direction problem, see Section 5 for details. Due to possible confounding, however, their method lacks a precise theoretical justification (Shimizu et al., 2006). We will see that it performs well only in low dimensions.

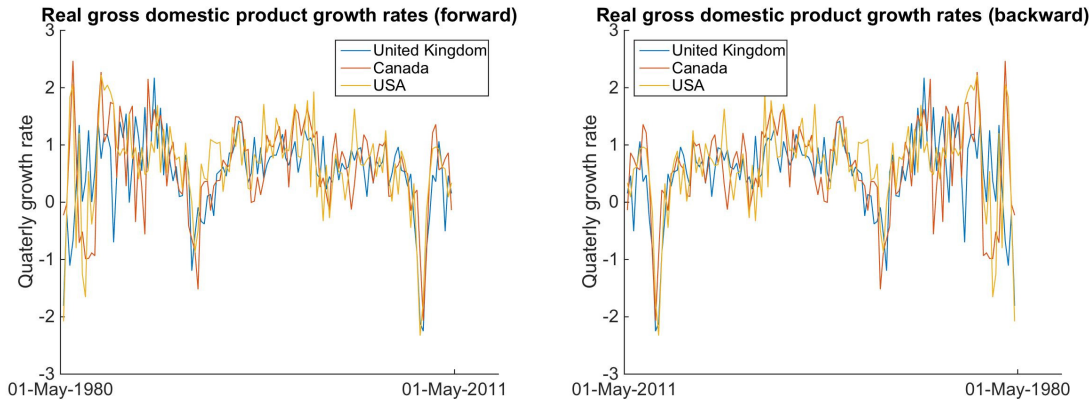

Figure 1: Quarterly growth rates in percentage of real gross domestic products for the United Kingdom, Canada and USA ((Tsay, 2014)). While the correct time direction is on the left, the time series on the right is reversed. The goal of this work is to investigate under which conditions the correct time direction becomes identifiable, without using any additional prior knowledge of the data domain.

The main idea of our method follows the univariate procedure of Peters et al. (2009). We consider vector-valued autoregressive moving average models (VARMA) with innovations that are independent of preceding values of the time series. We show that the time direction is reversible, i.e., the time series follows such a VARMA model in the reversed direction, only if the innovations are Gaussian. This leads to the following simple practical procedure for the determination of the arrow of time: after fitting a VARMA model to the observed time series, we test whether the residuals are independent of the preceding values in both time directions. If the model assumption is correct and the innovations are indeed non-Gaussian, we will find the independence only in one of the two possible directions.

In the remainder of Section 1, we introduce multivariate autoregressive moving average (VARMA) processes and formalize the problem. Section 2 contains identifiability results about VARMA processes of order one that are generalized to higher order processes in Section 3. We present a practical method in Section 4 and show results on simulated and real data sets in Section 5.

### 1.1 Notation and first properties of VAR models

A $K$-dimensional time series $\mathbf{X}_t$ is called a *vector autoregressive moving average process* VARMA$(p, q)$ if it is weakly stationary and if there is an i.i.d. white noise sequence $\mathbf{Z}_t$ with zero mean such that

$$\mathbf{X}_t - \mathbf{\Phi}_1 \mathbf{X}_{t-1} - \cdots - \mathbf{\Phi}_p \mathbf{X}_{t-p} = \mathbf{Z}_t + \mathbf{\Theta}_1 \mathbf{Z}_{t-1} + \cdots + \mathbf{\Theta}_q \mathbf{Z}_{t-q}. \tag{1}$$

Here, $\mathbf{X}_t$ being weakly stationary means that its mean is constant in time and the covariance between $\mathbf{X}_t$ and $\mathbf{X}_{t+h}$ depends only on the time lag $h \in \mathbb{Z}$. In more compact form, the above equation can be written as

$$\mathbf{\Phi}(B)\mathbf{X}_t = \mathbf{\Theta}(B)\mathbf{Z}_t, \tag{2}$$

where $\mathbf{\Phi}(z) = \mathbb{1} - \mathbf{\Phi}_1 z - \cdots - \mathbf{\Phi}_p z^p$, $\mathbf{\Theta}(z) = \mathbb{1} + \mathbf{\Theta}_1 z + \cdots + \mathbf{\Theta}_q z^q$ and $B$ is the backward shift operator. We call this process a VAR process for $p = 0$, and an MA process for $q = 0$. A VARMA$(p, q)$ process is further said to be *causal* if there exists a sequence of matrices $\mathbf{\Psi}_0, \mathbf{\Psi}_1, \ldots$ with absolutely summable components such that

$$\mathbf{X}_t = \sum_{j=0}^{\infty} \mathbf{\Psi}_j \mathbf{Z}_{t-j}. \tag{3}$$

The important condition here is that the index $j$ starts at zero, an expansion with $j \in \mathbb{Z}$ usually exists (Brockwell and Davis, 1991, Section 11.3). A way to check if a time series is causal is given by the following sufficient criterion:

**Lemma 1 (Causality criterion, Theorem 11.3.1 in (Brockwell and Davis, 1991))** *If $\det \mathbf{\Phi}(z) \neq 0 \; \forall z \in \mathbb{C}$ such that $|z| \leq 1$, then $\mathbf{\Phi}(B)\mathbf{X}_t = \mathbf{\Theta}(B)\mathbf{Z}_t$ with $\mathbf{Z}_t$ being white noise, has exactly one stationary solution, $\mathbf{X}_t = \sum_{j=0}^{\infty} \mathbf{\Psi}_j \mathbf{Z}_{t-j}$, where the matrices $\{\mathbf{\Psi}_j\}_{j \geq 0}$ have absolutely summable components and are uniquely determined by*

$$\mathbf{\Psi}(z) = \sum_{j=0}^{\infty} \mathbf{\Psi}_j z^j = \mathbf{\Phi}^{-1}(z)\mathbf{\Theta}(z), \quad |z| \leq 1.$$

In this work, we require yet another characterization of causal time series. As for all other results, the proof can be found in Section 6.

**Lemma 2** *A VARMA process is causal if and only if for all $i < t$ the noise $\mathbf{Z}_t$ is independent of the preceding value of the time series $\mathbf{X}_i$ (written as $\mathbf{Z}_t \perp\!\!\!\perp \mathbf{X}_i$).*

In this work, we investigate the time-reversibility of time series that follow such a causal VARMA model. More precisely, we call a causal VARMA process *time-reversible* if there is an i.i.d. noise sequence $\widetilde{\mathbf{Z}}_t$ and coefficient matrices $\widetilde{\boldsymbol{\Phi}}_1, \ldots, \widetilde{\boldsymbol{\Phi}}_p$ and $\widetilde{\boldsymbol{\Theta}}_1, \ldots, \widetilde{\boldsymbol{\Theta}}_p$ such that

$$\mathbf{X}_t - \widetilde{\boldsymbol{\Phi}}_1 \mathbf{X}_{t+1} - \cdots - \widetilde{\boldsymbol{\Phi}}_p \mathbf{X}_{t+p} = \widetilde{\mathbf{Z}}_t + \widetilde{\boldsymbol{\Theta}}_1 \widetilde{\mathbf{Z}}_{t+1} + \cdots + \widetilde{\boldsymbol{\Theta}}_q \widetilde{\mathbf{Z}}_{t+q}$$

where for all $i > t$, $\widetilde{\mathbf{Z}}_t \perp\!\!\!\perp \mathbf{X}_i$, see Lemma 2. In Sections 2 and 3, we try to find suitable criteria which allow us to determine whether a process is time-reversible. There, the independence between noise and values of the time series will play a crucial role.

### 1.2 Relation to the univariate case

Some special cases of VARMA processes relate directly to the univariate case. Some of these cases are degenerate and will be excluded from our analysis. In the case of a VAR(1) process with diagonal coefficient matrix of full rank, each component can be considered as a separate univariate time series and one can infer the true temporal ordering by directly applying the methods of Peters et al. (2009). For a two-dimensional VAR(1) process with coefficient matrix of the form $\boldsymbol{\Phi}_1 = \begin{bmatrix} \alpha & 0 \\ 0 & 0 \end{bmatrix}$, one can again apply the univariate results for the first element. Since the second time series component corresponds to i.i.d. noise, it is clear that this component is time-reversible independently of the noise distribution. In our analysis of the multivariate case, we will exclude the latter (degenerate) case, e.g., by assuming full rank of the coefficient matrix.

### 2. Time reversibility of VAR processes of order one

In the following subsection we only consider VAR(1) processes $\mathbf{X}_t = \boldsymbol{\Phi} \mathbf{X}_{t-1} + \mathbf{Z}_t$ and assume that $\boldsymbol{\Gamma}_0 := \operatorname{cov}(\mathbf{X}_t, \mathbf{X}_t)$ has full rank, see Section 1.2. This section contains the key theoretical argument of this work. Theorem 7 extends this argument to VARMA processes of any arbitrary but finite order.

We first show that we cannot infer the time direction of linear Gaussian time series. The proof is based on matrix algebra and the Yule-Walker equations.

**Proposition 3 (Gaussian errors lead to time-reversibility)** *Assume that the errors of a causal VAR(1) process $\mathbf{X}_t = \boldsymbol{\Phi} \mathbf{X}_{t-1} + \mathbf{Z}_t$ are normally distributed and that $\boldsymbol{\Gamma}_0 := \operatorname{cov}(\mathbf{X}_t, \mathbf{X}_t)$ is of full rank. Then, the process is time-reversible: there is a matrix $\widetilde{\boldsymbol{\Phi}}$ and a noise sequence $\widetilde{\mathbf{Z}}_t$ such that*

$$\mathbf{X}_t = \widetilde{\boldsymbol{\Phi}} \mathbf{X}_{t+1} + \widetilde{\mathbf{Z}}_t,$$

*where for all $i > t$, $\widetilde{\mathbf{Z}}_t$ is independent of $\mathbf{X}_i$.*

The following result is positive in the sense that non-Gaussian innovations introduce an asymmetry in the time direction.

**Theorem 4 (Non-Gaussian errors lead to time-identifiability)** *Consider a causal* VAR(1) *process $\mathbf{X}_t = \boldsymbol{\Phi} \mathbf{X}_{t-1} + \mathbf{Z}_t$ and assume that the process is time-reversible, i.e., there is a matrix $\widetilde{\boldsymbol{\Phi}}$ and a noise sequence $\widetilde{\mathbf{Z}}_t$ such that $\mathbf{X}_t$ can be written as:*

$$\mathbf{X}_t = \widetilde{\boldsymbol{\Phi}} \mathbf{X}_{t+1} + \widetilde{\mathbf{Z}}_t,$$

where for all $i > t$, $\widetilde{\mathbf{Z}}_t$ is independent of $\mathbf{X}_i$.

(i) *If the coefficient matrix $\mathbf{\Phi}$ is invertible, then all elements of the noise sequence vectors are normally distributed.*

(ii) *If the coefficient matrix $\mathbf{\Phi}$ is not nilpotent, then some elements of the noise sequence vectors are normally distributed.*

Our identifiability results come close to necessary and sufficient conditions. The case when $\mathbf{\Phi}$ is singular and not nilpotent and some but not all innovations are non-Gaussian is not covered. Both proofs can be found in Section 6. The Darmois-Skitovic theorem is used for proving the identifiability of independent component analysis (Comon, 1994) and lies at the heart of the proof of Theorem 4.

**Lemma 5 (Vectorised Darmois-Skitovich theorem for infinite sums)** *Let $\mathbf{X}_1, \mathbf{X}_2, \ldots$ be independent $d$-dimensional random vectors and consider the linear combinations $L_1 = \sum_{j=1}^{\infty} \mathbf{A}_j \mathbf{X}_j$ and $L_2 = \sum_{j=1}^{\infty} \mathbf{B}_j \mathbf{X}_j$ where $\mathbf{A}_j, \mathbf{B}_j$ are non-singular $d \times d$ matrices. If $L_1$ and $L_2$ are independent and $\{\mathbf{A}_j \mathbf{B}_j^{-1}\}_{j \geq 1}$ as well as $\{\mathbf{B}_j \mathbf{A}_j^{-1}\}_{j \geq 1}$ are bounded in some matrix norm, then the random vectors $\mathbf{X}_1, \mathbf{X}_2, \ldots$ are normally distributed.*

## 3. Time reversibility of VARMA processes of higher order

### 3.1 Definitions and parametrization

To ensure a unique solution to (1) in the univariate case, one requires that $\Phi(z)$ and $\Theta(z)$ of the univariate ARMA$(p,q)$ process $\Phi(B)X_t = \Theta(B)\epsilon_t$ have no common zeros (Bellmann, 1987). In order to guarantee a unique solution in the multivariate case, the conditions become slightly more complex.

As an example, consider a two-dimensional VARMA$(1,1)$ process

$$\mathbf{X}_t = \mathbf{\Phi}_1 \mathbf{X}_{t-1} + \mathbf{\Theta}_1 \mathbf{Z}_{t-1} + \mathbf{Z}_t$$

with $\mathbf{\Phi}_1 = \begin{bmatrix} 0 & \alpha + m \\ 0 & 0 \end{bmatrix}$ and $\mathbf{\Theta}_1 = \begin{bmatrix} 0 & -m \\ 0 & 0 \end{bmatrix}$, where $\alpha \neq 0$ and $m \in \mathbb{R}$. The same model can be written as a pure moving average process with $\mathbf{\Phi}_1 = \mathbf{0}$ and $\mathbf{\Theta}_1 = \begin{bmatrix} 0 & \alpha \\ 0 & 0 \end{bmatrix}$. Since $m$ in the VARMA$(1,1)$ presentation is arbitrary, the model parameters are not identifiable. This non-identifiability problem introduces one more difficulty to the problem of time-reversibility. In the univariate setting (Peters et al., 2009), one of the key ideas is to represent a VAR process as a MA$(\infty)$ process. For multivariate time series, however, the order of the corresponding MA process may be finite. For a VAR(1) process with $\mathbf{\Phi}_1 = \begin{bmatrix} 0 & 0.5 \\ 0 & 0 \end{bmatrix}$, for example, we find that $\mathbf{\Psi}_j = \mathbf{\Phi}^j = 0$ for $j > 1$ and thus $\mathbf{X}_t$ can be written as $\mathbf{X}_t = \mathbf{Z}_t + \mathbf{\Phi} \mathbf{Z}_{t-1}$, which is a pure MA process of finite order. These examples are taken from (Lütkepohl, 2010, Section 12.1.1).

Here, we solve this problem by assuming that $\mathbf{\Phi}_1$ is not nilpotent. Alternatively, one can require that the VARMA process is in the so called *final equations* or *echelon* form (Lütkepohl, 2010, Section 12.1.2).

## 3.2 Representing a VARMA process as a VAR process of order one

It is well known (e.g. Lütkepohl, 2010, Section 11.3.2) that a $K$-dimensional VARMA$(p,q)$ process

$$\mathbf{X}_t - \boldsymbol{\Phi}_1 \mathbf{X}_{t-1} - \cdots - \boldsymbol{\Phi}_p \mathbf{X}_{t-p} = \mathbf{Z}_t + \boldsymbol{\Theta}_1 \mathbf{Z}_{t-1} + \cdots + \boldsymbol{\Theta}_q \mathbf{Z}_{t-q} \qquad (4)$$

can be written as a $K(p+q)$-dimensional VAR(1) process

$$\widehat{\mathbf{X}}_t = \boldsymbol{\Upsilon} \widehat{\mathbf{X}}_{t-1} + \mathbf{U}_t,$$

with new noise innovations $\mathbf{U}_t$ and coefficients $\boldsymbol{\Upsilon}$ given by $\widehat{\mathbf{X}}_t = \begin{bmatrix} \mathbf{X}_t \ldots \mathbf{X}_{t-p+1} \mathbf{Z}_t \ldots \mathbf{Z}_{t-q+1} \end{bmatrix}^T$, $\mathbf{U}_t = \begin{bmatrix} \mathbf{Z}_t\, 0 \cdots 0\, \mathbf{Z}_t\, 0 \cdots 0 \end{bmatrix}^T$ each of dimension $K(p+q) \times 1$, and $\boldsymbol{\Upsilon} = \begin{bmatrix} \boldsymbol{\Upsilon}_{11} & \boldsymbol{\Upsilon}_{12} \\ \mathbf{0}_{Kq \times Kp} & \boldsymbol{\Upsilon}_{22} \end{bmatrix}$ with

$$\boldsymbol{\Upsilon}_{11} := \begin{bmatrix} \boldsymbol{\Phi}_1 & \boldsymbol{\Phi}_2 & \ldots & \boldsymbol{\Phi}_{p-1} & \boldsymbol{\Phi}_p \\ \mathbb{1}_K & 0 & \ldots & 0 & 0 \\ 0 & \mathbb{1}_K & \ldots & 0 & 0 \\ \vdots & & \ddots & 0 & \vdots \\ 0 & \ldots & 0 & \mathbb{1}_K & 0 \end{bmatrix}, \quad \boldsymbol{\Upsilon}_{12} := \begin{bmatrix} \boldsymbol{\Theta}_1 & \boldsymbol{\Theta}_2 & \ldots & \boldsymbol{\Theta}_{q-1} & \boldsymbol{\Theta}_q \\ 0 & 0 & \ldots & 0 & 0 \\ 0 & 0 & \ldots & 0 & 0 \\ \vdots & \vdots & & \vdots & \vdots \\ 0 & 0 & \ldots & 0 & 0 \end{bmatrix}$$

$$\boldsymbol{\Upsilon}_{22} := \begin{bmatrix} 0 & 0 & \ldots & 0 & 0 \\ \mathbb{1}_K & 0 & \ldots & 0 & 0 \\ 0 & \mathbb{1}_K & \ldots & 0 & 0 \\ \vdots & & \ddots & 0 & \vdots \\ 0 & \ldots & 0 & \mathbb{1}_K & 0 \end{bmatrix}.$$

The important property is that the first rows of $\boldsymbol{\Upsilon}$ equal $\begin{bmatrix} \boldsymbol{\Phi}_1 \cdots \boldsymbol{\Phi}_p\, \boldsymbol{\Theta}_1 \cdots \boldsymbol{\Theta}_q \end{bmatrix}$. With this and Lemma 1 we therefore have

**Lemma 6** *The VARMA$(p,q)$ process is causal if and only if its corresponding VAR(1) representation is causal.*

Thus the case of a VARMA$(p,q)$ process reduces to a VAR(1) process and we can apply the results from Section 2. Summarizing, we have the following theorem.

**Theorem 7** *Consider a VARMA$(p,q)$ process with not-nilpotent coefficient matrix $\boldsymbol{\Upsilon}$. If the error vectors are normally distributed, the process is time-reversible and if the process is time-reversible, then at least some of the elements of the noise vectors are normally distributed. (Note that $\boldsymbol{\Upsilon}$ is not nilpotent if and only if $\boldsymbol{\Phi}_1$ in the representation (4) is not nilpotent.)*

## 4. Algorithm

We now present a practical method for finite time series data and provide a consistency result for a version of this algorithm. For practical purposes, we restrict ourselves to VAR$(p)$ processes; see Lütkepohl (2010) for the technical difficulties with fitting VAR$(p,q)$ processes.

To estimate the correct direction of multivariate time series we follow Peters et al. (2009). The general idea is that under the discussed assumptions Gaussian causal VARMA processes are time-reversible but for any other error distribution we are able to identify the true temporal ordering, see Theorems 4 and 7. The main idea now is to fit VAR models in both directions and check in which direction we find that the residuals are independent of preceding values. The independence test is based on the Hilbert-Schmidt Independence Criterion (HSIC) (Gretton et al., 2007). In order to check whether the residual time series $\mathbf{Z}_t$ is independent of the past of $\mathbf{X}_t$, we simply check for independence between $\mathbf{Z}_t$ and $\mathbf{X}_{t-1}$ (although higher lags may be considered). The method decides correctly if the hypothesis of independence is rejected in the backward direction, while it is not rejected in the forward direction. In the case of Gaussian innovations, for example, we expect to accept the null hypotheses in both directions and the method remains undecided. For practical purposes, we introduce a small gap between the significance levels and include the option to work with the statistics itself rather with the p-value. The precise procedure is presented in Algorithm 1. For simulating (function *vgxproc*) and fitting (function *vgxvarx*) VAR($p$) processes we used

---

**Algorithm 1** Detecting the direction of multivariate time series

1: **procedure** INPUT($f = (x_1, \ldots, x_n), b = (x_n, \ldots, x_1)$, $sig1$, $sig2$),
2:    (a) fct := $-$HSIC, (b) fct := p $-$ valueHSIC
3:    $model_f = \text{VAR.fit}(f); \quad res_f = \text{residuals}(model_f); \quad fw = \text{fct}(f, res_f)$
4:    $model_b = \text{VAR.fit}(b); \quad res_b = \text{residuals}(model_b); \quad bw = \text{fct}(b, res_b)$
5:    **if** $\max(fw, bw) > sig1$ && $\min(fw, bw) < sig2$ **then**
6:       $decision = \text{argmax}(fw, bw)$
7:    **else** $decision =$ "I do not know."
8:    **end if**
9:    Return $decision$
10: **end procedure**

---

the "Econometrics Toolbox" within Matlab. The correct order of the process is estimated using AIC. The code is available at `http://people.inf.ethz.ch/bauers/`.

In practice, we might find that the model assumptions are violated due to the existence of hidden confounders or nonlinear relationships. Then we expect that the independence will be rejected in both directions and the method remains undecided. This is different from the non-decision in the Gaussian case, where we expect both directions to lead to a good model fit (see also Peters et al., 2009).

It can be shown that Algorithm 1 is consistent in the sense of Theorem 8 below. This result is not immediately straightforward since the independence measure is based on estimated residuals rather than "true" innovations (see also Mooij et al., 2014). Again, the proof is provided in Section 6.

**Theorem 8** *Let $(\mathbf{X}_t)_{t \in \mathbb{Z}}$ be a VAR process of order one with noise variables $(\mathbf{Z}_t)_{t \in \mathbb{Z}}$. Let the $\left(\mathbf{Z}_t^{fw}\right)_{t \in \mathbb{Z}}$ be the residuals in forward and $\left(\mathbf{Z}_t^{bw}\right)_{t \in \mathbb{Z}}$ the residuals in backward direction (corresponding to the best VAR fit) and assume that $\mathbf{X}_t \not\perp \mathbf{Z}_{t+1}^{bw}$ (see Theorem 4(i) and Remark 9 below). Assume that all processes are strictly stationary and uniformly mixing* with $\alpha(m) \leq m^{-3}$, as defined in (8). Then Algorithm 1 consistently estimates the arrow of time using an empirical HSIC score with a Gaussian kernel.

**Remark 9**
- For simplicity we assume that $\mathbf{X}_t \not\!\perp\!\!\!\perp \mathbf{Z}_{t+1}^{bw}$. One can use multiple testing to correct for a dependence at a different time lag.
- Under some technical assumptions Markov Chains (and thus ARMA processes) are uniformly mixing (e.g Doukhan, 1994; Mokkadem, 1988).
- The uniformly mixing assumption can be replaced by assuming that the process is absolutely regular (Chwialkowski and Gretton, 2014, Lemma 2).

## 5. Experiments

We compare our method with a LiNGAM-based approach proposed by Shimizu et al. (2006), which constructs a causal graph given i.i.d. samples of a random vector. If the generated graph is time consistent in the sense that all links go from lower to higher labeled variables (or the opposite), LiNGAM proposes this direction as the correct one.

### 5.1 Simulated data

We simulate VAR processes of different order and dimensionality and test the performance of both approaches. We use version (b) of our Algorithm 1 in order to better interpret our results. For all experiments, we use significance levels of $sig1 = 0.1$ and $sig2 = 0.05$. This is a conservative but interpretable choice that could be changed in order to increase the performance in the simulations.

**Simulation parameters** For a fixed parameter $\lambda = 2.5$, the $i$-th coefficient matrix of the simulated VAR process of dimension $k$ is generated by $\mathbf{\Phi}_i = \lambda^{-i}\mathbf{R} - (2\lambda)^{-i}\mathbf{Q}$, where $\mathbf{R}$ consists of uniformly drawn numbers between zero and one and $\mathbf{Q}$ is a matrix containing only ones.

**Deviation from Gaussian noise for different lags and dimensions** For $r$ ranging between zero to two we sampled each component of the noise vector as $Z_t \overset{\text{iid}}{\sim} \text{sgn}(Z) \cdot |Z|^r$, where $Z$ is Gaussian distributed. Only the case $r = 1$ results in a normal distribution and we should only then be unable to detect the true direction of the time series. This is verified for different lag orders $p$ and dimensions $k$ of the VAR process in Figure 2.

**Varying number of Gaussian error dimensions** Theorem 7 shows that one can detect the true direction of the time series if all error dimensions have a non-Gaussian distribution, while we can not infer the arrow if time when all errors are normally distributed. We therefore increased the number of Gaussian errors from 20% to 100%. Figure 3 supports our theoretical results and shows that our algorithm does not make a decision when all errors are Gaussian distributed. In addition, it suggests that with only one component not normally distributed we can still infer the true direction. This indicates that an even stronger version of Theorem 7 might hold.

**Comparison to LiNGAM** Figure 4 shows results of some of the settings shown in Figure 2. In general, the performance is comparable but the LiNGAM performance decreases for increasing dimension. Interestingly, LiNGAM makes no mistakes for $k = 3$ and $p = 3$ and Gaussian data. This might be due to the existence of $v$-structures in the graph (no instantaneous effects), see Figure 4(b). This is not the case for $k = 1$.

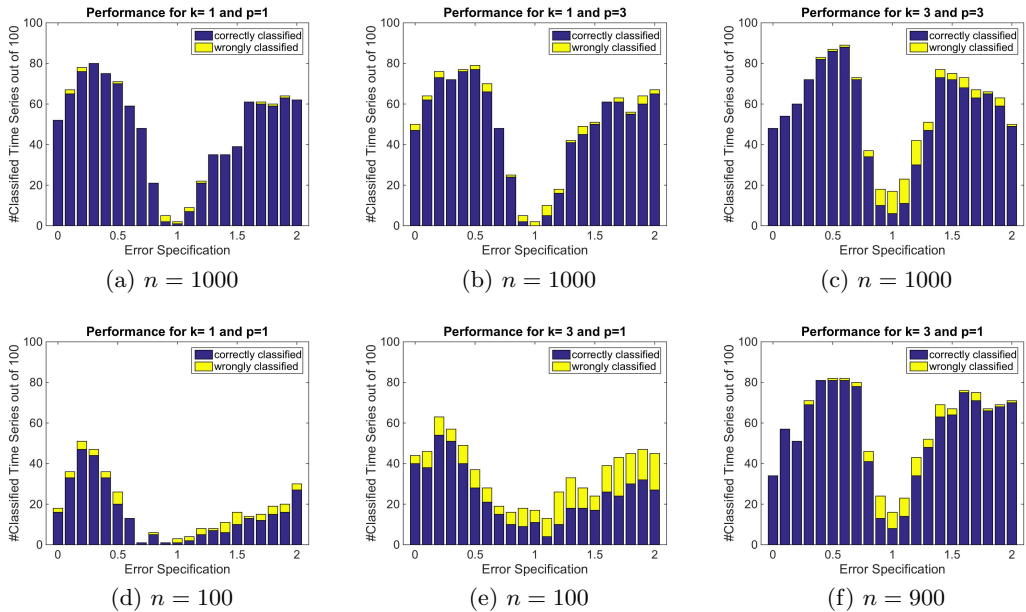

Figure 2: Correctly and incorrectly classified univariate time series for varying lag values, sample sizes and dimensions. For each value of $r$ we generated and tested our algorithm on 100 time series. Only $r = 1$ corresponds to Gaussian noise, for which the process is reversible. For these values our method remains undecided since both directions lead to a good model fit. It is interesting to note that the increased dimensionality ($k = 3$) introduces more model parameters but does not lead to worse performance with respect to identification of the time direction. If we compare situations, in which we have the same number of data points for each AR coefficient (namely 100, see second row (d) and (f)), the performance is significantly better for $k = 3$.

### 5.2 Real data: GDP growth for United Kingdom, Canada and United States

In a dataset containing the quarterly growth rates of real gross domestic product (GDP) of UK, Canada and USA from 1980 to 2011 (Tsay, 2014), we tested our approach for different time lags, see Figure 1. In Figure 5 the p-values of the independence test for forward and backward direction are plotted for time lags between one and ten. The optimal order chosen by AIC is four. If we treat the three time series individually, the method remains undecided in all cases. Only if we treat the process as a multivariate time series, the method outputs the correct result. The results do not change, when using version (a) of our Algorithm 1. LiNGAM does neither decide for the one-dimensional nor the three-dimensional case. Here, we took a first order difference which is often done in order to ensure stationarity.

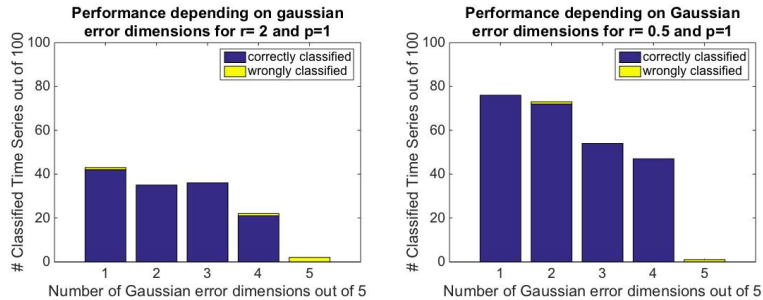

Figure 3: Dependence of our algorithm on the number of Gaussian error dimensions. A five dimensional VAR process is tested with increasing number of Gaussian errors. The remaining error element(s) follow a non-Gaussian noise distribution with $r = 0.5$ (cf. Section 5.1). In each case, the performance of the algorithm is tested on 100 time series with 1000 time points each. As expected, our method fails once the noise vector is completely normally distributed. As long as at least one element in the noise vector is non-Gaussian distributed, our method performs well.

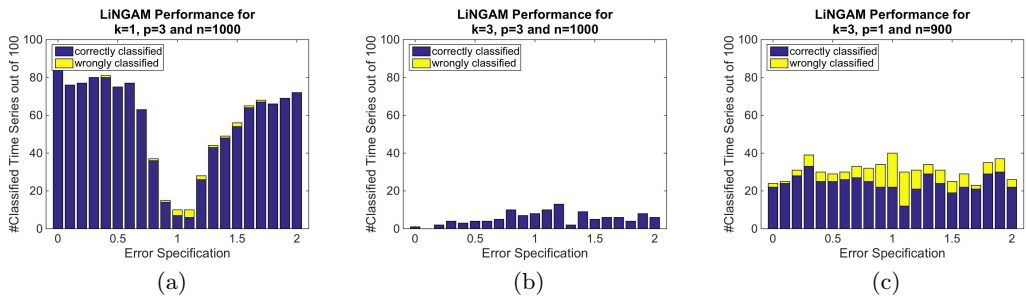

Figure 4: Correctly and wrongly classified time series with LiNGAM. While the performance seems to be slightly better in the univariate case, the performance decreases significantly with the dimension $k$ and lag order $p$ of the time series, compare with Figure 2.

### 5.3 Real data: video snippets Pickup et al. (2014)

In a computer vision application, Pickup et al. (2014) aim to determine if videos are shown in forward or backward direction. Apart from two other approaches (discriminative approach with training data and a heuristic claiming it is unlikely that multiple motions collapse into one motion) the authors apply an algorithm similar to the one presented by Peters et al. (2009) and our Algorithm 1(b). The authors model the velocities of moving points with a VAR(2) model and as outlined above, perform an independence test between velocities and model residuals (assuming that the univariate results of Peters et al. (2009) apply in higher dimensions). In correspondence to our results, the authors find that the approach works well if the assumptions like linear dynamics and non-Gaussian errors are satisfied.

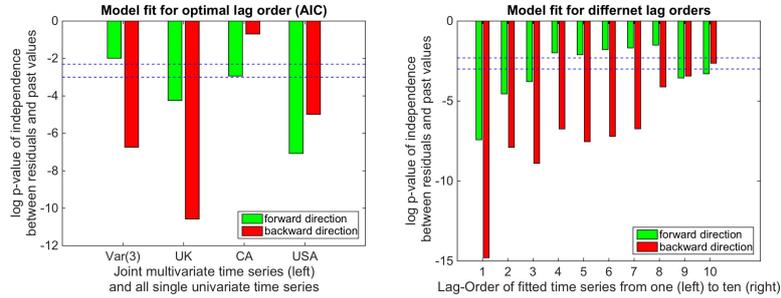

Figure 5: Left: log p-value for the independence test of residuals and time series values in forward and backward direction for the multivariate VAR process and the three univariate processes. For illustration purposes, dashed blue lines are shown at the values log(0.05) and log(0.1). Our method decides only if the p-values of both directions lie on different sides of this gap. While in the univariate case the algorithm does not take a decision, it is able to decides for the true temporal ordering in the multivariate case. The optimal orders are 4 for the multivariate process and $3, 5, 9$ in both directions for UK, Canada and USA. Right: the log p-values of the independence test for the multivariate process are shown for different orders. The performance of our algorithm is relatively robust to deviations from the optimal order.

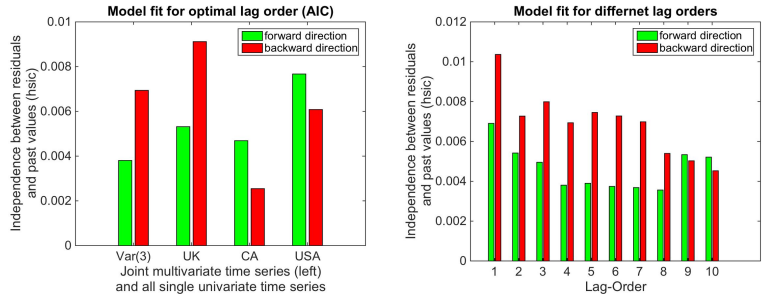

Figure 6: Left: HSIC values of the independence test of residuals and time series values in forward and backward direction for the multivariate VAR process and the three univariate processes. For the multivariate and the UK times series, the algorithm decides correctly while for the time series of Canada and the US, the algorithm decides wrongly. Right: the HSIC values of the independence test for the multivariate process are shown for different orders. Only for high lags and their inherent fitting problems the algorithm decides wrongly.

### 5.4 Real data: daily NASDAQ-stock values

In addition, we tested a set of 50 time series consisting of daily returns of NASDAQ stocks from the 1st October 2004 to 30th September 2007 (Yah, 2015). We grouped the time series into subsets of dimensions $k = 3, 5, 10$ and 50. For $k = 3$, $k = 5$ and $k = 10$ we chose not all but 100 randomly selected subsets. With $sig1 = 0.1$ and $sig2 = 0.05$ the method remained undecided in all cases. Apparently, the VAR model does not provide a good fit for these stock market data. This may be the case for other data sets, too, and we regard it as a positive feature that in those cases of model misspecification our method remains undecided rather than giving wrong answers.

The authors model the velocities of moving points with a VAR(2) model and as outlined above, perform an independence test between velocities and model residuals (assuming that the univariate results from Peters et al. (2009) also hold in higher dimensions). In correspondence to our results, the authors find that the approach works well if the assumptions like linearity of the dynamics and non-Gaussian errors are satisfied.

## 6. Proofs

### 6.1 Proof of Lemma 2

**Proof** The proof follows its univariate version (Peters et al., 2009). We first show that if $\mathbf{X}_t = \sum_{j=0}^{\infty} \mathbf{\Psi}_j \mathbf{Z}_{t-j}$, then $\mathbf{Z}_t$ is independent of $\mathbf{X}_i$, for all $i < t$. By defining $\mathbf{X}_t^{(n)} := \sum_{j=0}^{n} \mathbf{\Psi}_j \mathbf{Z}_{t-j}$, we have that $\mathbf{X}_t^{(n)}$ converges weakly to $\mathbf{X}_t$ and thus the characteristic function of $(\mathbf{X}_t, \mathbf{Z}_{t+1})$ obeys for all $u$ and $v$

$$\varphi_{\mathbf{P}(\mathbf{x}_t, \mathbf{z}_{t+1})}(u, v) = \lim_{n \to \infty} \varphi_{\mathbf{P}\left(\mathbf{x}_t^{(n)}, \mathbf{z}_{t+1}\right)}(u, v) = \lim_{n \to \infty} \varphi_{\mathbf{P}\mathbf{x}_t^{(n)}}(u) \varphi_{\mathbf{P}\mathbf{z}_{t+1}}(v) =$$
$$= \varphi_{\mathbf{P}\mathbf{x}_t}(u) \varphi_{\mathbf{P}\mathbf{z}_{t+1}}(v) = \varphi_{\mathbf{P}\mathbf{x}_t \otimes \mathbf{P}\mathbf{z}_{t+1}}(u, v).$$

This results in independence of $\mathbf{Z}_t$ and $X_i$, $i < t$ by uniqueness of the characteristic function.

In order to prove the "if"-part, we need to show that for a causal process $\mathbf{X}_t$, all coefficients $\mathbf{\Psi}_i$ are equal to zero for $i < 0$ in the Laurent expansion $\mathbf{X}_t = \sum_{i \in \mathbb{Z}} \mathbf{\Psi}_i \mathbf{Z}_{t-i}$ (see Lütkepohl (2010)). Assume otherwise, i.e., there is a coefficient $i_0 < 0$ such that $\mathbf{\Psi}_{i_0} \neq 0$. Then

$$\mathbf{\Psi}_{i_0} \mathbf{Z}_{t-i_0} + \sum_{i \in \mathbb{Z} - i_0} \mathbf{\Psi}_i \mathbf{Z}_{t-i} = \mathbf{X}_t \perp\!\!\!\perp \mathbf{\Psi}_{i_0} \mathbf{Z}_{t-i_0}. \tag{5}$$

Because $\mathbf{\Psi}_{i_0} \mathbf{Z}_{t-i_0}$ and $\sum_{i \in \mathbb{Z} - i_0} \mathbf{\Psi}_i \mathbf{Z}_{t-i}$ are independent with the same reasoning as above, (5) results in a contradiction. ∎

### 6.2 Proof of Proposition 3

**Proof** Defining

$$\widetilde{\mathbf{Z}}_t := \mathbf{X}_t - \text{cov}(\mathbf{X}_t, \mathbf{X}_{t+1}) \cdot \text{cov}(\mathbf{X}_{t+1}, \mathbf{X}_{t+1})^{-1} \mathbf{X}_{t+1}$$

and

$$\widetilde{\mathbf{\Phi}} := \text{cov}(\mathbf{X}_t, \mathbf{X}_{t+1}) \cdot \text{cov}(\mathbf{X}_{t+1}, \mathbf{X}_{t+1})^{-1},$$

it follows that

$$\widetilde{\mathbf{\Phi}} \mathbf{X}_{t+1} + \widetilde{\mathbf{Z}}_t = \text{cov}(\mathbf{X}_t, \mathbf{X}_{t+1}) \cdot \text{cov}(\mathbf{X}_{t+1}, \mathbf{X}_{t+1})^{-1} \mathbf{X}_{t+1}$$
$$+ \mathbf{X}_t - \text{cov}(\mathbf{X}_t, \mathbf{X}_{t+1}) \cdot \text{cov}(\mathbf{X}_{t+1}, \mathbf{X}_{t+1})^{-1} \mathbf{X}_{t+1} = \mathbf{X}_t.$$

In addition, we have

$$\begin{aligned}
\operatorname{cov}(\widetilde{\mathbf{Z}}_t, \mathbf{X}_{t+1}) &= \operatorname{cov}(\mathbf{X}_t - \operatorname{cov}(\mathbf{X}_t, \mathbf{X}_{t+1}) \cdot \operatorname{cov}(\mathbf{X}_{t+1}, \mathbf{X}_{t+1})^{-1} \mathbf{X}_{t+1}, \mathbf{X}_{t+1}) \\
&= \operatorname{cov}(\mathbf{X}_t, \mathbf{X}_{t+1}) - \operatorname{cov}(\mathbf{X}_t, \mathbf{X}_{t+1}) \cdot \operatorname{cov}(\mathbf{X}_{t+1}, \mathbf{X}_{t+1})^{-1} \operatorname{cov}(\mathbf{X}_{t+1}, \mathbf{X}_{t+1}) = 0.
\end{aligned}$$

By the assumption of the Gaussian distribution, the independence of $\widetilde{\mathbf{Z}}_t$ and $\mathbf{X}_{t+1}$ follows. It remains to show that $\widetilde{\mathbf{Z}}_t$ and $\mathbf{X}_{t+k}$ are independent for $k \geq 2$. By the mulitvariate form of the Yule-Walker equations for the VAR(1) process, i.e., $\boldsymbol{\Gamma}_k := \operatorname{cov}(\mathbf{X}_t, \mathbf{X}_{t+k}) = \boldsymbol{\Phi}\boldsymbol{\Gamma}_{k-1}$ (see 2.1.31 in Lütkepohl (2010)), we obtain that

$$\begin{aligned}
\operatorname{cov}(\widetilde{\mathbf{Z}}_t, \mathbf{X}_{t+k}) &= \operatorname{cov}(\mathbf{X}_t - \operatorname{cov}(\mathbf{X}_t, \mathbf{X}_{t+1}) \cdot \operatorname{cov}(\mathbf{X}_{t+1}, \mathbf{X}_{t+1})^{-1} \mathbf{X}_{t+1}, \mathbf{X}_{t+k}) \\
&= \operatorname{cov}(\mathbf{X}_t, \mathbf{X}_{t+k}) - \operatorname{cov}(\mathbf{X}_t, \mathbf{X}_{t+1}) \cdot \operatorname{cov}(\mathbf{X}_{t+1}, \mathbf{X}_{t+1})^{-1} \operatorname{cov}(\mathbf{X}_{t+1}, \mathbf{X}_{t+k}) \\
&= \operatorname{cov}(\mathbf{X}_t, \mathbf{X}_{t+k}) - \boldsymbol{\Phi}\boldsymbol{\Gamma}_0 \boldsymbol{\Gamma}_0^{-1} \operatorname{cov}(\mathbf{X}_{t+1}, \mathbf{X}_{t+k}) = \boldsymbol{\Gamma}_k - \boldsymbol{\Phi}\boldsymbol{\Gamma}_{k-1} = 0.
\end{aligned}$$

∎

### 6.3 Proof of Theorem 4

**Proof** We first prove the second part of the theorem, that is we assume that the $K \times K$ coefficient matrix $\boldsymbol{\Phi}$ is not nilpotent. By Lemma 1 it follows that

$$\mathbf{X}_t = \sum_{i=0}^{\infty} \boldsymbol{\Psi}_i \mathbf{Z}_{t-i} = \sum_{i=0}^{\infty} \boldsymbol{\Phi}^i \mathbf{Z}_{t-i} \qquad (6)$$

and

$$\widetilde{\mathbf{Z}}_{t-1} = \mathbf{X}_{t-1} - \widetilde{\boldsymbol{\Phi}} \mathbf{X}_t = \sum_{i=0}^{\infty} \boldsymbol{\Psi}_i \mathbf{Z}_{t-1-i} - \widetilde{\boldsymbol{\Phi}} \sum_{i=0}^{\infty} \boldsymbol{\Psi}_i \mathbf{Z}_{t-i} = \sum_{i=0}^{\infty} \left( \boldsymbol{\Psi}_{i-1} - \widetilde{\boldsymbol{\Phi}} \boldsymbol{\Psi}_i \right) \mathbf{Z}_{t-i}, \qquad (7)$$

where $\boldsymbol{\Psi}_{-1} := 0$, since for a VAR(1) process $\boldsymbol{\Psi}_i = \boldsymbol{\Phi}^i$. Defining $\mathbf{A}_i := \boldsymbol{\Psi}_{i-1} - \widetilde{\boldsymbol{\Phi}}\boldsymbol{\Psi}_i = (\mathbb{1} - \widetilde{\boldsymbol{\Phi}}\boldsymbol{\Phi})\boldsymbol{\Phi}^{i-1}$ and $\mathbf{B}_i := \boldsymbol{\Psi}_i = \boldsymbol{\Phi}^i = \boldsymbol{\Phi}\boldsymbol{\Phi}^{i-1}$, we can unfortunately not directly apply the vectorized Darmois-Skitovich theorem (see Theorem 4 of Ibragimov (2014) or our Lemma 5) since the matrices $\mathbf{A}_i$ and $\mathbf{B}_i$ are not invertible. By assumption the eigenvalues of the matrix $\boldsymbol{\Phi}$ of a causal process are smaller than one in absolute value (see Lemma 1) and thus by Gelfand's formula $(\mathbb{1} - \widetilde{\boldsymbol{\Phi}}\boldsymbol{\Phi})$ is of full rank $K$. Using the assumption that $\boldsymbol{\Phi}$ is not nilpotent it follows that for a large enough index $i \geq i_0$ and some number $s$, we have $\operatorname{rank}(\boldsymbol{\Phi}^i) = \operatorname{rank}(\boldsymbol{\Phi}^{i-1}) = s$. Applying the rank inequality of Sylvester, $(\operatorname{rank}(\mathbf{Q}) + \operatorname{rank}(\mathbf{P}) - K \leq \operatorname{rank}(\mathbf{QP}) \leq \min(\operatorname{rank}(\mathbf{Q}), \operatorname{rank}(\mathbf{P}))$ with $\mathbf{Q}$ and $\mathbf{P}$ matrices of size $K \times K$), we get that $\operatorname{rank}(\mathbf{A}_i) = s$ and $\operatorname{rank}(\mathbf{B}_i) = s$ for $i \geq i_0$.

A singular value decomposition yields $\boldsymbol{\Phi}^{i-1} = \mathbf{U}\boldsymbol{\Sigma}\mathbf{V}^*$ with a unitary matrix $\mathbf{U}$ and a $K \times K$ diagonal matrix $\boldsymbol{\Sigma}$, whose first $s$ diagonal elements are non-zero, while its last $K - s$ columns contain only zeros. We define a new noise vector $\boldsymbol{\varepsilon}_{t-i} := \mathbf{V}^* \mathbf{Z}_{t-i}$ and rewrite (6) and (7) as

$$\mathbf{X}_t = \sum_{i=0}^{\infty} \boldsymbol{\Psi}_i \mathbf{Z}_{t-i} = \sum_{i=0}^{\infty} \boldsymbol{\Phi}\mathbf{U}\boldsymbol{\Sigma}\mathbf{V}^* \mathbf{Z}_{t-i} = \sum_{i=0}^{\infty} \boldsymbol{\Phi}\mathbf{U}\boldsymbol{\Sigma} \boldsymbol{\varepsilon}_{t-i}$$

and

$$\widetilde{\mathbf{Z}}_{t-1} = \sum_{i=0}^{\infty} \left(\mathbb{1} - \widetilde{\boldsymbol{\Phi}}\boldsymbol{\Phi}\right) \mathbf{U}\boldsymbol{\Sigma}\mathbf{V}^*\mathbf{Z}_{t-i} = \sum_{i=0}^{\infty} \left(\mathbb{1} - \widetilde{\boldsymbol{\Phi}}\boldsymbol{\Phi}\right) \mathbf{U}\boldsymbol{\Sigma}\boldsymbol{\varepsilon}_{t-i}.$$

Since $\text{rank}(\boldsymbol{\Phi}^{i-1}) = s$, $\left(\mathbb{1} - \widetilde{\boldsymbol{\Phi}}\boldsymbol{\Phi}\right)\mathbf{U}\boldsymbol{\Sigma}$ and $\boldsymbol{\Phi}\mathbf{U}\boldsymbol{\Sigma}$ are of the form $(\mathbf{M}_1|0)$ and $(\mathbf{M}_2|0)$ for $K \times s$ matrices $\mathbf{M}_j$ with $\text{rank}(\mathbf{M}_j) = s$, $j \in \{1,2\}$. Let us select $s$ independent rows of $\mathbf{M}_2$ and call them $\mathbf{B}_i^*$. Let us call $\mathbf{A}_i^*$ the corresponding elements in $\mathbf{M}_1$. This leads to a new, reduced system of equations

$$\mathbf{Z}_{t-1}^* = \sum_{i=0}^{\infty} \mathbf{A}_i^* \begin{pmatrix} \varepsilon^1 \\ \varepsilon^2 \\ \vdots \\ \varepsilon^s \end{pmatrix}_{t-i} \quad \text{and} \quad \mathbf{X}_t^* = \sum_{i=0}^{\infty} \mathbf{B}_i^* \begin{pmatrix} \varepsilon^1 \\ \varepsilon^2 \\ \vdots \\ \varepsilon^s \end{pmatrix}_{t-i},$$

where $\mathbf{A}_i^*$ and $\mathbf{B}_i^*$ are of dimension $s \times s$ with full rank, and $\mathbf{Z}_{t-1}^*$ and $\mathbf{X}_t^*$ are independent. Darmois-Skitovich (Lemma 5 in the supplement) implies normality of $\varepsilon_t^1, \ldots, \varepsilon_t^s$ and the normality of some elements of the original noise vector $\mathbf{Z}_t$ follows by Cramer's theorem.

For the second part of the theorem let us now assume that the coefficient matrix $\boldsymbol{\Phi}$ and thus $\mathbf{A}_i$ and $\mathbf{B}_i$ are invertible (by the rank inequality of Sylvester). Since

$$\mathbf{A}_i\mathbf{B}_i^{-1} = (\mathbb{1} - \widetilde{\boldsymbol{\Phi}}\boldsymbol{\Phi})\boldsymbol{\Phi}^{i-1}(\boldsymbol{\Phi}^i)^{-1} = (\mathbb{1} - \widetilde{\boldsymbol{\Phi}}\boldsymbol{\Phi})\boldsymbol{\Phi}^{-1} = \boldsymbol{\Phi}^{-1} - \widetilde{\boldsymbol{\Phi}}$$

and

$$\mathbf{B}_i\mathbf{A}_i^{-1} = (\mathbf{A}_i\mathbf{B}_i^{-1})^{-1} = (\boldsymbol{\Phi}^{-1} - \widetilde{\boldsymbol{\Phi}})^{-1}$$

do not depend on $i$, both $\{\mathbf{A}_i\mathbf{B}_i^{-1}\}_{i \geq 1}$ and $\{\mathbf{B}_i\mathbf{A}_i^{-1}\}_{i \geq 1}$ are bounded w.r.t. some norm. The Darmois-Skitovich theorem (Lemma 5 in the supplement) then implies the normality of $\mathbf{Z}_t$. ∎

### 6.4 Proof of Lemma 6

**Proof** The statement follows from considering determinants, see Lemma 1.

$$\det(\mathbb{1}_{K(p+q)} - \boldsymbol{\Upsilon}z) = \det(\mathbb{1}_{Kp} - \boldsymbol{\Upsilon}_{11}z)\det(\mathbb{1}_{Kq} - \boldsymbol{\Upsilon}_{22}z) = \det(\mathbb{1}_K - \boldsymbol{\Phi}_1 z - \cdots - \boldsymbol{\Phi}_p z^p)$$

∎

### 6.5 Proof of Theorem 8

**Proof** Our proof follows the main argument of Theorem 20 in Mooij et al. (2014) and omits some of the details. We first report a consistency results for HSIC for time dependent data (Chwialkowski and Gretton, 2014).

Given a stationary process $(W_t)_{t \in \mathbb{Z}} := (\mathbf{X}_t, \mathbf{Z}_t)_{t \in \mathbb{N}}$, let $(W_t^\star)_{t \in \mathbb{Z}}$ be a sequence of independent copies of $W_0$. Let $\mathbf{W}_{\underline{T}} := (W_1, \ldots, W_T)$ for a sample of size $T$. We show that the empirical HSIC (with a predefined, data independent bandwidth) between estimated residuals and the time series converges to its true value under some mixing condition. A process is uniformly mixing with $\alpha(m)$ if

$$\alpha(m) := \sup_n \sup_{A \in \mathrm{A}_1^n} \sup_{B \in \mathrm{A}_{n+m}^\infty} |P(B|A) - P(B)| \to 0, \tag{8}$$

where $\mathrm{A}_b^c = \sigma(W_b, W_{b+1}, \ldots, W_c)$ is a sigma field spanned by $W_b, W_{b+1}, \ldots, W_c$ (Chwialkowski and Gretton, 2014).

For $w_j := (x_j, z_j)$, $j \in \{a, b, c, d\}$ and $k$ and $\ell$ being Gaussian kernels, define

$$h(w_a, w_b, w_c, w_d) := k(x_a, x_b)[\ell(z_a, z_b) + \ell(z_c, z_d) - 2\ell(z_b, z_c)].$$

We further define

$$\gamma := \mathrm{HSIC}_{\mathbf{X}_0, \mathbf{Z}_0} = \mathrm{E}[h(W_1^\star, W_2^\star, W_3^\star, W_4^\star)].$$

Since Gaussian kernels are characteristic, it is known that $\mathbf{X}_0$ is independent of $\mathbf{Z}_0$ if and only if $\gamma = 0$ (Gretton et al., 2007). We therefore have for any $t \in \{1, \ldots, T\}$,

$$\mathbf{X}_t \perp\!\!\!\perp \mathbf{Z}_t \iff \mathbf{X}_0 \perp\!\!\!\perp \mathbf{Z}_0 \iff \mathrm{E}[h(W_1^\star, W_2^\star, W_3^\star, W_4^\star)] = 0.$$

An empirical estimate $\widehat{\mathrm{HSIC}}(\mathbf{W}_{\underline{T}})$ of $\gamma$ based on the sample $\mathbf{W}_{\underline{T}}$ can be calculated through the $V$-statistic

$$\widehat{\mathrm{HSIC}}(\mathbf{W}_{\underline{T}}) := V_T(h, \mathbf{W}_{\underline{T}}) := \frac{1}{T^4} \sum_{1 \leq t_1, t_2, t_3, t_4 \leq T} h(W_{t_1}, W_{t_2}, W_{t_3}, W_{t_4}).$$

Note that unlike in the *i.i.d.* case, here the $W_t$ are not independent for different values of $t$. However, we still have the following consistency result (requiring the uniformly mixing assumption in this step): For $T \to \infty$ we have for independent $\mathbf{X}_t$ and $\mathbf{Z}_t$

$$V_T(h, \mathbf{W}_{\underline{T}}) \xrightarrow{P} 0 \tag{9}$$

and for $\mathbf{X}_t$ and $\mathbf{Z}_t$ being dependent:

$$V_T(h, \mathbf{W}_{\underline{T}}) \xrightarrow{P} \delta > 0. \tag{10}$$

This follows from Theorems 1 and 2 in Chwialkowski and Gretton (2014) using that for any random variable $\mathbf{S}$ with finite variance $T \cdot \mathbf{X}_T \xrightarrow{D} \mathbf{S}$ implies that $\mathbf{X}_T \xrightarrow{L^2} 0$ and therefore $\mathbf{X}_T \xrightarrow{P} 0$ for $T \to \infty$.

Let now $\widehat{\mathbf{Z}}_{\underline{T}}^{\mathrm{fw}} := (\widehat{\mathbf{Z}}_2^{\mathrm{fw}}, \ldots, \widehat{\mathbf{Z}}_{T+1}^{\mathrm{fw}})$ be the residuals in forward and $\widehat{\mathbf{Z}}_{\underline{T}}^{\mathrm{bw}} := (\widehat{\mathbf{Z}}_2^{\mathrm{bw}}, \ldots, \widehat{\mathbf{Z}}_{T+1}^{\mathrm{bw}})$ the residuals in backward direction. We will show that for $\mathbf{X}_{\underline{T}} = (\mathbf{X}_1, \ldots, \mathbf{X}_T)$,

$$\widehat{\mathrm{HSIC}}\left(\mathbf{X}_{\underline{T}}, \widehat{\mathbf{Z}}_{\underline{T}}^{\mathrm{fw}}\right) \xrightarrow{P} 0 \text{ and } \widehat{\mathrm{HSIC}}\left(\mathbf{X}_{\underline{T}}, \widehat{\mathbf{Z}}_{\underline{T}}^{\mathrm{bw}}\right) \xrightarrow{P} \delta > 0. \tag{11}$$

By (a slightly modified version of) Lemma 16 in Mooij et al. (2014) it follows that

$$\left|\widehat{\mathrm{HSIC}}(\mathbf{X}_{\underline{T}}, \widehat{\mathbf{Z}}_{\underline{T}}) - \widehat{\mathrm{HSIC}}(\mathbf{X}_{\underline{T}}, \mathbf{Z}_{\underline{T}})\right|^2 \leq \left(\frac{C}{\sqrt{T}}\right)^2 \|\widehat{\mathbf{Z}}_{\underline{T}} - \mathbf{Z}_{\underline{T}}\|_F^2, \tag{12}$$

for some constant $C \in \mathbb{R}$ and Frobenius norm $\|\cdot\|_F$.

Given the "correct" VAR(1) representation $\mathbf{X}_t = \mathbf{\Phi}\mathbf{X}_{t-1} + \mathbf{Z}_t$ and the model fit $\mathbf{X}_t = \widehat{\mathbf{\Phi}}\mathbf{X}_{t-1} + \widehat{\mathbf{Z}}_t$, we have

$$\mathrm{E}\left[\frac{1}{T}\|\widehat{\mathbf{Z}}_{\underline{T}} - \mathbf{Z}_{\underline{T}}\|_F^2\right] = \mathrm{E}\left[\frac{1}{T}\sum_t \|\widehat{\mathbf{\Phi}}\mathbf{X}_{t-1} - \mathbf{\Phi}\mathbf{X}_{t-1}\|^2\right] \leq \mathrm{E}\left[\|\widehat{\mathbf{\Phi}} - \mathbf{\Phi}\|_F^2 \frac{1}{T}\sum_t \|\mathbf{X}_{t-1}\|^2\right]$$

Since the variance of $\mathbf{Z}_t$ and thus $\mathbf{X}_{t-1}$ is assumed to be finite and parameter estimation in VAR processes is consistent (Lütkepohl, 2010, Chapter 3), the above expression goes to zero for $T \to \infty$. Since the right hand side in (12) vanishes asymptotically in expectation, it follows that

$$\lim_{T \to \infty} \mathrm{E}\left[\left|\widehat{\mathrm{HSIC}}(\mathbf{X}_{\underline{T}}, \widehat{\mathbf{Z}}_{\underline{T}}) - \widehat{\mathrm{HSIC}}(\mathbf{X}_{\underline{T}}, \mathbf{Z}_{\underline{T}})\right|^2\right] = 0.$$

Since convergence in $L_2$ implies convergence in probability it follows that

$$\widehat{\mathrm{HSIC}}(\mathbf{X}_{\underline{T}}, \widehat{\mathbf{Z}}_{\underline{T}}) - \widehat{\mathrm{HSIC}}(\mathbf{X}_{\underline{T}}, \mathbf{Z}_{\underline{T}}) \xrightarrow{P} 0.$$

Together with

$$\widehat{\mathrm{HSIC}}(\mathbf{X}_{\underline{T}}, \mathbf{Z}_{\underline{T}}) \xrightarrow{P} \mathrm{HSIC}(\mathbf{X}_0, \mathbf{Z}_0)$$

this implies

$$\widehat{\mathrm{HSIC}}(\mathbf{X}_{\underline{T}}, \widehat{\mathbf{Z}}_{\underline{T}}) \xrightarrow{P} \mathrm{HSIC}(\mathbf{X}_0, \mathbf{Z}_0) \tag{13}$$

Since $\mathbf{X}_t \not\perp\!\!\!\perp \mathbf{Z}^{\mathrm{bw}}_{t+1}$ our statement (11) follows by combining the above convergence result (13) with (9) and (10). ∎

## 7. Discussion and future work

We have derived a framework for the identification of the direction of multivariate time series. By assuming that the data generating process exhibits linear dynamics and non-Gaussian noise we were able to extend the results of Peters et al. (2009) to multiple dimensions. In addition, we provide a consistency result for our algorithm that covers those of Peters et al. (2009) and Pickup et al. (2014) as special cases. The approach works well on simulated and some financial data. The empirical results for simulated data sets indicate that the detection of the temporal ordering of a time sequence might be possible as long as one element of the noise vector is normally distributed, which would be slightly stronger than the theoretical guarantee we provide. In general, the performance for determining the direction of real world time series (e.g., video snippets, see Pickup et al. (2014)) depends on the validity of our assumptions; in particular, this includes linear dynamics and non-Gaussian additive noise. We found that in the case of model misspecification, our approach usually remains undecided rather than giving incorrect answers. An extension of our framework to non-linear dynamics could reduce the number of non-decisions.